\documentclass{amsart}
\usepackage{amsmath, amscd, amsthm}
\usepackage{latexsym, amssymb, amsfonts}
\usepackage{url}
\usepackage{color}
\usepackage{graphicx}
\usepackage{hyperref}
\hypersetup{
    colorlinks=false,
    linkcolor=blue,
    filecolor=magenta,      
    urlcolor=cyan,
    pdftitle={Fundamental domain for the Markoff-Hurwitz equation},
    pdfpagemode=FullScreen,
}

\usepackage{stmaryrd}
\usepackage{geometry}\usepackage{tikz}\usetikzlibrary{positioning}
\usepackage{float}
\usepackage{subcaption}

\numberwithin{equation}{section}

\def\Bigroman{\uppercase\expandafter{\romannumeral\number\count 255 }}
\def\Romannumeral{\afterassignment\BigRoman\count255=}

\theoremstyle{definition}
\newtheorem{definition}{Definition}[section]

\theoremstyle{remark}
\newtheorem*{remark}{Remark}

\theoremstyle{plain}
\newtheorem{theorem}{Theorem}[section] 

\newtheorem{lemma}[theorem]{Lemma}
\newtheorem{corollary}[theorem]{Corollary}

\title{Fundamental domain for the Markoff-Hurwitz equation}
\author{Eunju Shin}
\address{Department of Mathematics, Seoul National University, 1 Gwanak-ro, Gwanak-gu, Seoul 08826}
\email{eshin@snu.ac.kr}
\usepackage{lipsum}

\begin{document}
\begin{abstract} For integers $a\neq0$, $k$, and $n\geq3$, we consider the Markoff-Hurwitz equation given by $x_1^1+\cdots+x_n^2-ax_1\cdots x_n=k$. By defining graphs associated with a height function and using their properties, we find an exact fundamental domain for a symmetric group generated by involution maps sending $(x_1,\dots,x_n)$ to $(x_1,\dots,ax_1\cdots x_{i-1}x_{i+1}\cdots x_n-x_i,\dots,x_n)$, permutations, and double sign changes on the set of integral solutions for the Markoff-Hurwitz equation.
    
\end{abstract}
\maketitle

\section{Introduction} 

Historically, the integral solutions for the Markoff equation given by $x^2+y^2+z^2-xyz=0$ have been one of the classical problems in Diophantine equations. Ghosh and Sarnak's  recent study\cite{gs} is about a fundamental set for the symmetric group generated by the Vieta involution maps $(x,y,z)\mapsto(yz-x,y,z)$, $(x,y,z)\mapsto(x,zx-y,z)$, and $(x,y,z)\mapsto(z,y,xy-z)$, permutations, and double sign changes on the set of integral solutions for the Markoff equation $x^2+y^2+z^2-xyz=k$ for some $k\neq4$. One way to generalize the Markoff equation, namely the Markoff-Hurwitz equation, is by increasing the number of variables in the Markoff equation and multiplying its specific term by a constant.

For $a,k,n\in\mathbb{Z}$ with $a\neq0$ and $n\geq3$, the \textit{Markoff-Hurwitz equation} is given by
\begin{align*}
    V_{a,k,n}: x_1^2+\cdots+x_n^2-ax_1\cdots x_n=k,
\end{align*}
and denote by $V_{a,k,n}(\mathbb{Z})$ the set of its integral solutions. For $i=1,\dots,n$, define involution maps 
\begin{align*}
    \mathcal{V}_{a,n,i}(x_1,\dots,x_n)=(x_1,\dots,ax_1\cdots x_{i-1}x_{i+1}\cdots x_n-x_i,\dots,x_n).
\end{align*}
For $1\leq s<t\leq n$, let $\mathcal{S}_{n,s,t}(x_1,\dots,x_n)=(x_1,\dots,-x_s,\dots,-x_t,\dots,x_n)$ (double sign change). For $\sigma\in S_n$, let $\sigma(x_1,\dots,x_n)=(x_{\sigma(1)},\dots,x_{\sigma(n)})$ (permutation). Let 
\begin{align*}
    \Gamma_{a,n}=\langle\mathcal{V}_{a,n,i},\mathcal{S}_{n,s,t},\sigma:i,s,t\in\{1,
    \dots,n\}, s<t, \sigma\in S_n\rangle,
\end{align*}
 which acts on $V_{a,k,n}(\mathbb{Z})$ ($\mathbf{x}\mapsto \gamma\cdot\mathbf{x}=\gamma(\mathbf{x})$ for $\gamma\in\Gamma_{a,n}$) because the equation is symmetric,
 \begin{align*}
 x_1^2+\cdots+x_{i-1}^2+x_{i+1}^2+\cdots+x_n^2-x_i(ax_1\cdots x_{i-1}x_{i+1}\cdots x_n-x_i)&=k,\ \text{and}\\
 x_1^2+\cdots+(-x_s)^2+\cdots+(-x_t)^2+\cdots+x_n^2-ax_1\cdots(-x_s)\cdots(-x_t)\cdots x_n&=k.
 \end{align*}
From the equation $(-x_1)^2+x_2^2+\cdots-(-a)(-x_1)x_2\cdots x_n=x_1^2+\cdots-ax_1\cdots x_n$, by replacing $x_1$ with $-x_1$, it suffices to treat when $a>0$. The main result of this paper is as follows:
 \begin{theorem}\label{1.1} Let $a$, $k$, and $n$ be integers such that $a>0$ and $n\geq3$. The set
     \begin{align*}
    \mathfrak{S}_0(a,k,n)\cup\mathfrak{S}_1(a,k,n)\cup\mathfrak{S}_2(a,k,n)\cup\mathfrak{S}_{>2}(a,k,n),
\end{align*}
where 
\begin{align*}
    \mathfrak{S}_0(a,k,n)&=\{(0,x_2,\dots,x_n)\in\mathbb{Z}^n:0\leq x_2\leq\cdots\leq x_n, x_2^2+\cdots+x_n^2=k\},\\
    \mathfrak{S}_1(a,k,n)&=\left\{\begin{matrix}\{(-1,1,\dots,1,x_{n-1},x_n)\in \mathbb{Z}^n:1\leq x_{n-1}\leq x_n,
 x_{n-1}^2+x_n^2+x_{n-1}x_n=k-n+2\}&\text{if}\ a=1\\
\emptyset\ \ \ \ \ \ \ \ \ \ \ \ \ \ \ \ \ \ \ \ \ \ \ \ \ \ \ \ \ \ \ \ \ \ \ \ \ \ \ \ \ \ \ \ \ \ \ \ \ \ \ \ \ \ \ \ \ \ \ \ \ \ \ \ \ \ \ \ \ \ \ \ \ \ \ \ \ \ \ \ \ \ \ \ \ \ \ \ \ \ \ \ \ \ \ \ \ & \text{if}\ a>1\end{matrix}\right.,\\
    \mathfrak{S}_2(a,k,n)&=\left\{\begin{matrix}\begin{matrix}\{(1,\dots,1,2,x_{n-1},x_{n-1})\in\mathbb{Z}^n:2\leq x_{n-1},k=n+1\}\ \ \ \ \ \ \ \ \ \ \ \ \ \ \ \ \ \ \ \ \ \ \ \ \ \ \ \ \ \ \\
        \cup\{(-1,1,\dots,1,2,x_{n-1},x_n)\in \mathbb{Z}^n:2\leq x_{n-1}\leq x_n, x_{n-1}+x_n=\sqrt{k-n-1}\}\end{matrix} & \text{if}\ a=1\\
\begin{matrix}\{(1,\dots,1,x_{n-1},x_{n-1})\in\mathbb{Z}^n:1\leq x_{n-1},k=n-2\}\ \ \ \ \ \ \ \ \ \ \ \ \ \ \ \ \ \ \ \ \ \ \ \ \ \ \ \ \ \ \\
        \cup\{(-1,1,\dots,1,x_{n-1},x_n)\in \mathbb{Z}^n:1\leq x_{n-1}\leq x_n, x_{n-1}+x_n=\sqrt{k-n+2}\}\end{matrix} \ \ \ & \text{if}\ a=2\\
\emptyset\ \ \ \ \ \ \ \ \ \ \ \ \ \ \ \ \ \ \ \ \ \ \ \ \ \ \ \ \ \ \ \ \ \ \ \ \ \ \ \ \ \ \ \ \ \ \ \ \ \ \ \ \ \ \ \ \ \ \ \ \ \ \ \ \ \ \ \ \ \ \ \ \ \ \ \ \ \ \ \ \ \ \ \ \ \ \ \ \ \ \ \ & \text{if}\ a>2
\end{matrix}\right.
,\\
    \mathfrak{S}_{>2}(a,k,n)&=\left\{\mathbf{x}\in V_{k,n}(\mathbb{Z}):0<x_1\leq\cdots\leq x_n\leq\frac{a}{2}x_1\cdots x_{n-1}\ \text{and}\ ax_1\cdots x_{n-2}>2\right\}\\
        &\ \ \ \ \cup \{\mathbf{x}\in V_{k,n}(\mathbb{Z}):0<-x_1\leq x_2\leq \cdots\leq x_n\ \text{and}\ ax_1\cdots x_{n-2}<-2\},
\end{align*}
is a fundamental domain for $\Gamma_{a,n}$ on $V_{a,k,n}(\mathbb{Z})$.
 \end{theorem}
As its corollary, we conclude that for $(a,k)\neq (1,n+1), (2,n-2)$ and $a>0$, $V_{a,k,n}(\mathbb{Z})$ has at most finitely many $\Gamma_{a,n}$-orbits. We will use a method similar to \cite{gs} and \cite{s} for the proof, introducing graphs associated with a height function and observing their properties.

\subsection*{Acknowledgements} I would like to thank my advisor, Junho Peter Whang, and my colleague, Seongyoon Kim, for introducing me to this problem and for their valuable contributions during our helpful discussions.


\section{Fundamental domain for the Markoff-Hurwitz equation}
Let $a,k,n\in\mathbb{Z}$ with $a>0$ and $n\geq3$. Let $V_{a,k,n}, V_{a,k,n}(\mathbb{Z}),\mathcal{V}_{a,n,i},\mathcal{S}_{n,s,t}$, and $\Gamma_{a,n}$ be the same as above. 
\subsection{Graphs associated with a height function} 
Define a height function 
\begin{align*}
    \Delta(x_1,\dots,x_n)=|x_1|+\cdots+|x_n|.
\end{align*}
Now, define a graph $\mathcal{G}$ as follows: each integral solution $\mathbf{x}\in V_{a,k,n}(\mathbb{Z})/\sim$, where $\mathbf{x}\sim\mathcal{S}_{n,s,t}(\mathbf{x})$ for some $1\leq s<t\leq n$ or  $\mathbf{x}\sim\sigma(\mathbf{x})$ for some $\sigma\in S_n$, will be a vertex. The vertices are distributed in the plane in order of the values of $\Delta$ (vertices with smallest values are located at the bottom) and each vertex $\mathbf{x}$ is connected with vertices $\mathcal{V}_{a,n,i}(\mathbf{x})$ by edges. 
\begin{remark} For any $\mathbf{x}\in V_{a,k,n}(\mathbb{Z})$, $1\leq s<t\leq n$, and $i\neq s,t$, we have
    \begin{align*}
        \mathcal{V}_{a,n,i}\mathcal{S}_{n,s,t}(\mathbf{x})&=(x_1,\dots,-x_s,\dots,-x_t,\dots,ax_1\cdots(-x_s)\cdots(-x_t)\cdots x_{i-1}x_{i+1}\cdots x_n-x_i,\dots,x_n)\\
        &=(x_1,\dots,-x_s,\dots,-x_t,\dots,ax_1\cdots x_{i-1}x_{i+1}\cdots x_n-x_i,\dots,x_n)\\
        &=\mathcal{S}_{n,s,t}(x_1,\dots,ax_1\cdots x_{i-1}x_{i+1}\cdots x_n-x_i,\dots,x_n)=\mathcal{S}_{n,s,t}\mathcal{V}_{a,n,i}(\mathbf{x}),\\
        \mathcal{V}_{a,n,s}\mathcal{S}_{n,s,t}(\mathbf{x})&=(x_1,\dots,ax_1\cdots x_{s-1}x_{s+1}\cdots(-x_t)\cdots x_n-(-x_s),\dots,-x_t,\dots,x_n)\\
        &=(x_1,\dots,-(ax_1\cdots x_{s-1}x_{s+1}\cdots x_n-x_s),\dots,-x_t,\dots,x_n)\\
        &=\mathcal{S}_{n,s,t}(x_1,\dots,ax_1\cdots x_{s-1}x_{s+1}\cdots x_n-x_s,\dots,x_n)=\mathcal{S}_{n,s,t}\mathcal{V}_{a,n,s}(\mathbf{x}),
    \end{align*}
    and similarly, $\mathcal{V}_{a,n,t}\mathcal{S}_{n,s,t}(\mathbf{x})=\mathcal{S}_{n,s,t}\mathcal{V}_{a,n,t}(\mathbf{x})$. For any $\mathbf{x}\in V_{a,k,n}(\mathbb{Z})$ and $\sigma\in S_n$, we have
    \begin{align*}
        \mathcal{V}_{a,n,i}\sigma(\mathbf{x})&=(x_{\sigma(1)},\dots,x_{\sigma(i-1)},ax_{\sigma(1)}\cdots x_{\sigma(i-1)}x_{\sigma(i+1)}\cdots x_{\sigma(n)}-x_{\sigma(i)},x_{\sigma(i+1)},\dots,x_{\sigma(n)})\\
        &=\sigma(x_1,\dots,x_{\sigma(i)-1},ax_{\sigma(1)}\cdots x_{\sigma(i-1)}x_{\sigma(i+1)}\cdots x_{\sigma(n)}-x_{\sigma(i)},x_{\sigma(i)+1},\dots,x_n)=\sigma\mathcal{V}_{a,n,\sigma(i)}(\mathbf{x}).
    \end{align*}
Consequently (reordering edges if necessary), we conclude that if $\mathbf{x}\sim\mathbf{y}$, then  $\mathcal{V}_{a,n,i}(\mathbf{x})\sim\mathcal{V}_{a,n,i}(\mathbf{y})$ for all $i=1,\dots,n$. Hence, the above graphs are well-defined. 
\end{remark}

\begin{definition}\label{def:2.1}  Let $\mathcal{G}$ be a graph associated with $\Delta$. Let $\mathbf{x}$ and $\mathbf{y}$ be vertices on $\mathcal{G}$. 
\begin{itemize}
    \item If there are $l_1,\dots,l_N\in\{1,\dots,n\}$ such that $\mathbf{y}=\mathcal{V}_{a,n,l_N}\cdots\mathcal{V}_{a,n,l_1}(\mathbf{x})$ and  
\begin{align*}
    \Delta\mathcal{V}_{a,n,l_{t+1}}\mathcal{V}_{a,n,l_t}\cdots\mathcal{V}_{a,n,l_1}(\mathbf{x})\geq\Delta\mathcal{V}_{a,n,l_t}\cdots\mathcal{V}_{a,n,l_1}(\mathbf{x})\geq\Delta(\mathbf{x})
\end{align*}
for all $i\in\{1,\dots,n\}$ and $t\in\{1,\dots,N-1\}$, then we say $\mathbf{x}$ is the \textit{descendent vertex} of $\mathbf{y}$.  
\item If  $\Delta(\mathbf{x})<\Delta\mathcal{V}_{a,n,i}(\mathbf{x})$ for all $i=1,\dots,n$, then we say $\mathbf{x}$ is the \textit{last vertex} of $\mathcal{G}$.  
\end{itemize}
\end{definition}

\subsection{Properties of the graphs associated with the height function}

Consider the set
\begin{align*}
    \mathfrak{T}(a,k,n)=\{\mathbf{x}\in V_{a,k,n}(\mathbb{Z}):|ax_{i_1}\cdots x_{i_{n-2}}|\leq2\ \text{for some}\ 1\leq i_1<\cdots<i_{n-2}\leq n\}.
\end{align*}
\begin{lemma}\label{2.1}
    Let $\mathbf{x}\in V_{a,k,n}(\mathbb{Z})\setminus\mathfrak{T}(a,k,n)$. If $|ax_1\cdots x_{i-1}x_{i+1}\cdots x_n-x_i|\leq |x_i|$, then $|x_j|<|x_j|$ for all $j\neq i$.
\end{lemma}
\textit{Proof}. Since $|ax_1\cdots x_{i-1}x_{i+1}\cdots x_n-x_i|\leq |x_i|$, then by the triangular inequality, we have $|ax_1\cdots x_{i-1}x_{i+1}\cdots x_n|\leq2|x_i|$. Since $|ax_{i_1}\cdots x_{i_{n-2}}|>2$ holds for all $1\leq i_1<\cdots<i_{n-2}\leq n$, it follows that for all $j\neq i$,
\begin{align*}
2|x_j|&<|ax_1\cdots x_{j-1}x_{j+1}\cdots x_{i-1}x_{i+1}\cdots x_n||x_j|=|ax_1\cdots x_j\cdots x_{i-1}x_{i+1}\cdots x_n|\leq 2|x_i|,
\end{align*}
that is, $|x_j|<|x_i|$.\qed

\begin{remark}
    If $\mathbf{x}\in V_{a,k,n}(\mathbb{Z})\setminus\mathfrak{T}(a,k,n)$, then $ax_1\cdots x_{i-1}x_{i+1}\cdots x_n-x_i\neq-x_i$ for all $i=1,\dots,n$ because, on the contrary, if $ax_1\cdots x_{i-1}x_{i+1}\cdots x_n-x_i=-x_i$, then for some $j\neq i$,
\begin{align*}
    2<|ax_1\cdots x_{j-1}x_{j+1}\cdots x_{i-1}x_{i+1}\cdots x_n|\leq|ax_1\cdots x_{i-1}x_{i+1}\cdots x_n|=0,
\end{align*}
which is a contradiction. Then by Lemma \ref{2.1}, for $\mathbf{x}\in V_{a,k,n}(\mathbb{Z})\setminus\mathfrak{T}(a,k,n)$, there are two cases:\footnotesize
\begin{align*}
    \begin{tikzpicture}[every node/.style={align=center}] 
    \node(a){$\mathbf{x}$}; 
    \node(b)[below =of a]{$\mathcal{V}_{a,n,i}(\mathbf{x})$} edge [-](a);
    \node(c)[above= of a]{$\cdots$} edge [-] (a);
    \node(d)[above left= of a]{$\mathcal{V}_{a,n,1}(\mathbf{x})$} edge [-] (a);
    \node(e)[above right= of a]{$\mathcal{V}_{a,n,n}(\mathbf{x})$} edge [-] (a);
    \end{tikzpicture}\hspace{8pt}
    \begin{tikzpicture}[every node/.style={align=center}] 
    \node(a){$\mathbf{x}$}; 
    \node(b)[above= of a]{$\cdots$} edge [-] (a);
    \node(c)[above left= of a]{$\mathcal{V}_{a,n,1}(\mathbf{x})$} edge [-] (a);
    \node(d)[above right= of a]{$\mathcal{V}_{a,n,n}(\mathbf{x})$} edge [-] (a);
    \end{tikzpicture}
\end{align*}\small
\end{remark}
Note that $\Delta$ is invariant under double sign change and permutation, so we may assume that $0\leq|x_1|\leq x_2\leq\cdots\leq x_n$, and replace $\mathfrak{T}(a,k,n)$ with
\begin{align*}
    \mathfrak{T}'(a,k,n)&=\{\mathbf{x}\in V_{a,k,n}(\mathbb{Z}):0\leq|x_1|\leq x_2\leq\cdots\leq x_n\ \text{and}\ |ax_1\cdots x_{n-2}|\leq2\}\\
    &=\mathfrak{T}_0'(a,k,n)\cup\mathfrak{T}_1'(a,k,n)\cup\mathfrak{T}_2'(a,k,n)
\end{align*}
where
\begin{align*}
    \mathfrak{T}_s'(a,k,n)=\{\mathbf{x}\in\mathfrak{T}'(a,k,n):|ax_1\cdots x_{n-2}|=s\}\ \text{for}\ s=0,1,2.
\end{align*} 
Note that 
\begin{align*}
\mathfrak{T}_1'(a,k,n)&=\left\{\mathbf{x}\in\mathfrak{T}'(a,k,n):|x_1\cdots x_{n-2}|=\frac{1}{a}\right\}=\emptyset\ \text{if}\ a>1,\\
\mathfrak{T}_2'(a,k,n)&=\left\{\mathbf{x}\in\mathfrak{T}'(a,k,n):|x_1\cdots x_{n-2}|=\frac{2}{a}\right\}=\emptyset\ \text{if}\ a>2.
\end{align*}


\begin{lemma}\label{2.2} Let $\mathbf{x}\in\mathfrak{T}_1'(1,k,n)$ with $x_1=1$. Then either \begin{enumerate}
    \item[(\romannumeral1)] $x_{n-1}=x_n$ and $\mathbf{x}$ has a descendent vertex $\mathbf{y}=(0,1,\dots,1,x_{n-1})\in \mathfrak{T}_0'(1,k,n)$, or
    \item[(\romannumeral2)] $x_{n-1}<x_n$ and $\mathbf{x}$ has a descendent vertex $\mathbf{y}=(-1,1,\dots,1,x_n-x_{n-1},x_{n-1})\in\mathfrak{T}_1'(1,k,n)$.
\end{enumerate}
Moreover, $\mathbf{y}\sim\mathcal{V}_{1,n,n}(\mathbf{x})$ is a descendent vertex of $\mathbf{x},\mathcal{V}_{1,n,i}(\mathbf{x})$ for all $i=1,\dots,n$.
\end{lemma}
\textit{Proof}. Assume $x_1=1$. Then $x_1=\cdots=x_{n-2}=1\leq x_{n-1}\leq x_n$.
\begin{enumerate}
    \item[(\romannumeral1)] If $x_{n-1}=x_n$, then $x_1\cdots x_{n-2}x_n-x_{n-1}=x_{n-1}-x_n=0$, implying that 
    \begin{align*}
        \mathcal{V}_{1,n,n-1}(\mathbf{x})\sim\mathcal{V}_{1,n,n}(\mathbf{x})\sim(0,1,\dots,1,x_{n-1})\in \mathfrak{T}_0'(1,k,n).
    \end{align*} 
    For all $i\leq n-2$, we have $x_1\cdots x_{i-1}x_{i+1}\cdots x_n-x_i=x_{n-1}^2-1\geq0$. Since $x_{n-1}^2-1=0$ implies $x_{n-1}=1$, then $\mathcal{V}_{1,n,i}(\mathbf{x})\sim(0,1,\dots,1)=(0,1,\dots,1,x_{n-1})$, it follows that for some $j\leq n-2$, a part near $\mathbf{x}$ of the graph associated with $\Delta$ having $\mathbf{x}$ as a vertex is shown in the diagrams below:\footnotesize
\begin{figure}[H]
\centering
\subcaptionbox{\footnotesize$x_{n-1}>1$\normalsize}{\begin{tikzpicture}[every node/.style={align=center}] 
    \node(a){$\mathbf{x}=(0,1,\dots,x_{n-1},x_{n-1})$};
    \node(b)[below = of a]{$\begin{matrix}\mathcal{V}_{a,n,n-1}(\mathbf{x})\sim\mathcal{V}_{a,n,n-1}(\mathbf{x})\\
\sim(0,1,\dots,1,x_{n-1})\in\mathfrak{T}_0'(1,k,n)\end{matrix}$} edge [-] (a);
\node(c)[above = of a]{$\begin{matrix}\mathcal{V}_{a,n,1}(\mathbf{x})\sim\cdots\mathcal{V}_{a,n,n-2}(\mathbf{x})\\
\sim(0,1,\dots,1,x_{n-1},x_{n-1},x_{n-1}^2-1)\end{matrix}$} edge [-] (a);
    \end{tikzpicture}}%
\hspace{8pt}
\subcaptionbox{\footnotesize$x_{n-1}=1$\normalsize}{\begin{tikzpicture}[every node/.style={align=center}] 
    \node(a){$\mathbf{x}=(1,\dots,1)$};
    \node(b)[below = of a]{$\mathcal{V}_{a,n,i}(\mathbf{x})\sim(0,1,\dots,1)\in\mathfrak{T}_0'(1,k,n)$} edge [-] (a);
    \end{tikzpicture}}
\end{figure}\normalsize

    \item[(\romannumeral2)] If $x_{n-1}<x_n$, then $x_1\cdots x_{n-1}-x_n=x_{n-1}-x_n<0$ and 
    \begin{align*}
        \mathcal{V}_{1,n,n}\mathcal{V}_{1,n,n-1}(\mathbf{x})\sim\mathcal{V}_{1,n,n}(\mathbf{x})\sim(-1,1,\dots,1,x_n-x_{n-1},x_{n-1})\in\mathfrak{T}_1'(1,k,n).
    \end{align*}
     For all $i\leq n-2$, we have $x_1\cdots x_{i-1}x_{i+1}\cdots x_n-x_i=x_{n-1}x_n-1\geq1$.
     \begin{itemize}
         \item If $x_{n-1}x_n-1=1$, then $x_{n-1}=1$, $x_n=2$, $\mathcal{V}_{1,n,n-1}(\mathbf{x})=\mathcal{V}_{1,n,i}(\mathbf{x})=\mathbf{x}$.\footnotesize
         \begin{align*}
             \begin{tikzpicture}[every node/.style={align=center}] 
    \node(a){$\mathbf{x}=(1,\dots,1,1,2)$};
    \node(b)[below = of a]{$\mathcal{V}_{1,n,n}(\mathbf{x})\sim(-1,1,\dots,1,1,1)
    \in\mathfrak{T}_1'(1,k,n)$} edge [-] (a);
    \end{tikzpicture}
         \end{align*}\normalsize
         \item If $x_{n-1}x_n-1>1$, then $\Delta(\mathbf{x})<\Delta\mathcal{V}_{1,n,i}(\mathbf{x})$ for all $i\leq n-2$, and $x_n\geq3$. Then we have $(x_n-x_{n-1})x_n-1\geq2>1$, which means that $\Delta(\mathcal{V}_{1,n,n-1}(\mathbf{x}))<\Delta\mathcal{V}_{1,n,i}(\mathcal{V}_{1,n,n-1}(\mathbf{x}))$ for all $i\leq n-2$. Since $\Delta\mathcal{V}_{1,n,n}\mathcal{V}_{1,n,n-1}(\mathbf{x})<\Delta\mathcal{V}_{1,n,n-1}(\mathbf{x})<\Delta(\mathbf{x})$, then we have the following diagram:\footnotesize
         \begin{align*}
     \begin{tikzpicture}[every node/.style={align=center}] 
    \node(a){$\mathbf{x}=(1,\dots,1,x_{n-1},x_n)$}; 
    \node(b)[above= of a]{$\mathcal{V}_{1,n,1}(\mathbf{x})\sim\cdots\sim\mathcal{V}_{1,n,n-2}(\mathbf{x})$} edge [-] (a);
    \node(c)[below left= of a]{$\mathcal{V}_{1,n,n-1}(\mathbf{x})=(1,\dots,1,x_n-x_{n-1},x_n)$} edge [-] (a);
\node(f)[below = of a]{};
    \node(d)[below = of f]{$\begin{matrix}
        \mathcal{V}_{1,n,n}\mathcal{V}_{1,n,n-1}(\mathbf{x})
        \sim\mathcal{V}_{1,n,n}(\mathbf{x})\\
        \sim(-1,1,\dots,1,x_n-x_{n-1},x_{n-1})
    \in\mathfrak{T}_1'(1,k,n)
    \end{matrix}$} edge [-] (a) edge [-](c);
\node(e)[left= of a]{$\begin{matrix}\mathcal{V}_{1,n,1}\mathcal{V}_{1,n,n-1}(\mathbf{x})\sim\mathcal{V}_{1,n,n-2}\mathcal{V}_{1,n,n-1}(\mathbf{x})\\
\sim(1,\dots,1,x_n-x_{n-1},(x_n-x_{n-1})x_n-1,x_n)\end{matrix}
$} edge [-] (c);
    \end{tikzpicture}
\end{align*}\normalsize
     \end{itemize}
     Thus, $\mathcal{V}_{1,n,n}(\mathbf{x})$ is a descendent vertex of $\mathbf{x},\mathcal{V}_{1,n,i}(\mathbf{x})$ for all $i=1,\dots,n$.\qed
\end{enumerate}
\begin{remark}
    Let $\mathbf{x}\in\mathfrak{T}_1'(1,k,n)$ with $x_1=-1$. Then $-x_1=x_2=\cdots=x_{n-2}=1\leq x_{n-1}\leq x_n$ and for $1<i$,
    \begin{align*}
        |x_1\cdots x_{i-1}x_{i+1}\cdots x_n-x_i|&>|x_i|,\\
        |x_2\cdots x_n-x_1|=|x_{n-1}x_n-1|&>1=|x_1|,
    \end{align*}
    implying that $\Delta(\mathbf{x})<\Delta\mathcal{V}_{1,n,j}(\mathbf{x})$ for all $j=1,\dots,n$. That is, $\mathbf{x}$ is the last vertex of the graph associated with $\Delta$, corresponding to the orbit $\Gamma_{1,n}\cdot\mathbf{x}$. Then by Lemma \ref{2.2}, we conclude that $\mathbf{y}\in\mathfrak{T}_1'(1,k,n)$ is either 
    \begin{enumerate}
        \item[(\romannumeral1)] a descendent vertex of $\mathbf{y},\mathcal{V}_{1,n,i}(\mathbf{y})$ for all $i=1,\dots,n$, or 
        \item[(\romannumeral2)] the last vertex of the graph associated with $\Delta$, corresponding to the orbit $\Gamma_{1,n}\cdot\mathbf{y}$.
    \end{enumerate}
\end{remark}


\subsection{Proof of Theorem \ref{1.1}} Now, we collect all the last vertices of each $\Gamma_{a,n}$-orbit and claim that the set of these last vertices is a fundamental domain for $\Gamma_{a,n}$. For this, we need to show the following lemma.
\begin{lemma}\label{2.3} Let $\mathbf{x}\in V_{a,k,n}(\mathbb{Z})$. The orbit $\Gamma_{a,n}\cdot\mathbf{x}$ has exactly one last vertex.    
\end{lemma}
\textit{Proof}. Let $\mathbf{x}\in V_{a,k,n}(\mathbb{Z})$ and $\mathbf{y}\in\Gamma_{a,n}\cdot\mathbf{x}$ with $|y_1|\leq y_2\leq\cdots\leq y_n$. Let 
\begin{align*}
    \mathfrak{T}_{>2}'(a,k,n)=\{\mathbf{z}\in V_{a,k,n}(\mathbb{Z}):|z_1|\leq z_2\leq\cdots z_n\ \text{and}\ |az_1\cdots z_{n-2}|>2\}.
\end{align*}
\begin{enumerate}
\item[\romannumeral1)] If $\mathbf{y}\in\mathfrak{T}_0'(a,k,n)$, then for $i>1$,
\begin{align*}
\mathcal{V}_{a,n,1}(\mathbf{y})&=\mathcal{V}_{a,n,1}(0,y_2,\dots,y_n)=(y_2\cdots y_n,y_2,\dots,y_n),\\
\mathcal{V}_{a,n,1}(\mathbf{y})&=\mathcal{V}_{a,n,i}(0,y_2,\dots,y_n)=(0,y_2,\dots,-y_i,\dots,y_n).
\end{align*}
Thus, $\Delta\mathcal{V}_{a,n,1}(\mathbf{y})-\Delta(\mathbf{y})=y_2\cdots y_n,y_2\geq0$,  $\Delta\mathcal{V}_{a,n,i}(\mathbf{y})=\Delta(\mathbf{y})$ for $i>1$, and there are two cases, as shown in the following diagrams:\footnotesize
\begin{figure}[H]
\centering
\subcaptionbox{\footnotesize$y_1=y_2=0$\normalsize}{\begin{tikzpicture}[every node/.style={align=center}] 
    \node(a){$\mathbf{y}=(0,0,y_3,\dots,y_n)\sim\mathcal{V}_{a,n,i}(\mathbf{y})$};
    \end{tikzpicture}}
\hspace{8pt}
\subcaptionbox{\footnotesize$y_1=0$ and $y_2\neq0$\normalsize}{\begin{tikzpicture}[every node/.style={align=center}] 
    \node(a){$\mathbf{y}=(0,y_2,\dots,y_n)\sim\mathcal{V}_{a,n,i}(\mathbf{y})$};
    \node(b)[above = of a]{$\mathcal{V}_{a,n,1}(\mathbf{y})$} edge [-] (a);
    \end{tikzpicture}}
\end{figure}\normalsize
\item[\romannumeral2)] If $\mathbf{y}\in\mathfrak{T}_1'(a,k,n)\cup\mathfrak{T}_{>2}'(a,k,n)$, since $\mathfrak{T}_1'(a,k,n)=\emptyset$ for all $a>1$, by Lemma \ref{2.1}, \ref{2.2}, and the above remark, we have that $\mathbf{y}\in\mathfrak{T}_1'(a,k,n)\cup\mathfrak{T}_{>2}'(a,k,n)$ is either 
    \begin{enumerate}
        \item[(\romannumeral1)] a descendent vertex of $\mathbf{y},\mathcal{V}_{a,n,i}(\mathbf{y})$ for all $i=1,\dots,n$, or 
        \item[(\romannumeral2)] the last vertex of the graph associated with $\Delta$, corresponding to the orbit $\Gamma_{a,n}\cdot\mathbf{x}$.
    \end{enumerate}

\item[\romannumeral3)] If $\mathbf{y}\in\mathfrak{T}_2'(a,k,n)$, since $\mathfrak{T}_2'(a,k,n)=\emptyset$ for all $a>2$, we may assume that $a\leq2$.
\begin{enumerate}
\item[(\romannumeral1)] Assume $a=1$. It follows that $|y_1|=y_2=\cdots=y_{n-3}=1<y_{n-2}=2\leq y_{n-1}\leq y_n$, and we have tha for $j\leq n-3$, 
\begin{align*}
|y_1\cdots y_{j-1}y_{j+1}\cdots y_n-y_j|&\geq2|y_{n-1}y_n|-1>|y_j|,\\
|y_1\cdots y_{n-3}y_{n-1}y_n-y_{n-2}|&\geq |y_{n-1}y_n|-2\geq|y_{n-2}|,\\
|y_1\cdots y_{n-2}y_n-y_{n-1}|&\geq2|y_n|-|y_{n-1}|\geq|y_{n-1}|.
\end{align*}
Note that if $|y_1y_{n-1}y_n-2|\leq2$, since $2\leq y_{n-1}\leq y_n$ and $|y_1\cdots y_{n-3}y_{n-1}y_n|>2$, then $y_1=\cdots=y_{n-3}=1<y_{n-2}=y_{n-1}=y_n=2$.  It follows that
    \begin{align*}
        &|y_1\cdots y_{n-3}y_{n-1}y_n-y_{n-2}|>|y_{n-2}|\ \text{or}\\
        &y_1=\cdots=y_{n-3}=1<y_{n-2}=y_{n-1}=y_n=2.
    \end{align*}
    Similarly, if $|2y_1y_n-y_{n-1}|\leq y_{n-1}$, since $2\leq y_{n-1}\leq y_n$ and $|y_1\cdots y_{n-2}y_n|>2$, then $y_1=\cdots=y_{n-3}=1<y_{n-2}=2\leq y_{n-1}=y_n$. Thus, we have that
    \begin{align*}
        &|y_1\cdots y_{n-3}y_{n-2}y_n-y_{n-1}|>|y_{n-1}|\ \text{or}\\
        &y_1=\cdots=y_{n-3}=1<y_{n-2}=2\leq y_{n-1}=y_n.
    \end{align*}
    That is, there are four cases, as shown in the diagrams below:  \footnotesize
    \begin{figure}[H]
    \centering
    \subcaptionbox{\footnotesize$\Delta\mathcal{V}_{1,n,n}(\mathbf{y})<\Delta(\mathbf{y})<\Delta\mathcal{V}_{1,n,n-1}(\mathbf{y})$\normalsize}{\begin{tikzpicture}[every node/.style={align=center}] 
    \node(a){$\mathbf{y}$}; 
    \node(b)[below =of a]{$\mathcal{V}_{1,n,n}(\mathbf{y})$} edge [-](a);
    \node(c)[above= of a]{$\cdots$} edge [-] (a);
    \node(d)[above left= of a]{$\mathcal{V}_{1,n,1}(\mathbf{y})$} edge [-] (a);
    \node(e)[above right= of a]{$\mathcal{V}_{1,n,n-1}(\mathbf{y})$} edge [-] (a);
    \end{tikzpicture}}
\hspace{8pt} 
\subcaptionbox{\footnotesize$\Delta(\mathbf{y})<\Delta\mathcal{V}_{1,n,n-1}(\mathbf{y}),\Delta\mathcal{V}_{1,n,n}(\mathbf{y})$\normalsize}{\begin{tikzpicture}[every node/.style={align=center}] 
    \node(a){$\mathbf{y}$};
    \node(c)[above= of a]{$\cdots$} edge [-] (a);
    \node(d)[above left= of a]{$\mathcal{V}_{1,n,1}(\mathbf{y})$} edge [-] (a);
    \node(e)[above right= of a]{$\mathcal{V}_{1,n,n}(\mathbf{y})$} edge [-] (a);
    \end{tikzpicture}}
\vspace{10pt}

\subcaptionbox{\footnotesize$\Delta(\mathbf{y})=\Delta\mathcal{V}_{1,n,n-1}(\mathbf{y})<\Delta\mathcal{V}_{1,n,n-2}(\mathbf{y})$\normalsize}{\begin{tikzpicture}[every node/.style={align=center}] 
    \node(a){$\mathbf{y}=(1,\dots,1,2,x_{n-1},x_{n-1})$}; 
    \node(c)[above = of a]{$\mathcal{V}_{1,n,1}(\mathbf{y})\sim\cdots\sim\mathcal{V}_{1,n,n-3}(\mathbf{y})$} edge [-] (a);
    \node(d)[above right= of a]{$\mathcal{V}_{1,n,n-2}(\mathbf{y})$} edge [-] (a);
    \end{tikzpicture}}
\subcaptionbox{\footnotesize$\Delta(\mathbf{y})=\Delta\mathcal{V}_{1,n,n-1}(\mathbf{y})=\Delta\mathcal{V}_{1,n,n-2}(\mathbf{y})$\normalsize}{\begin{tikzpicture}[every node/.style={align=center}] 
    \node(a){$\mathbf{y}=(1,\dots,1,2,2,2)$}; 
    \node(b)[above = of a]{$\mathcal{V}_{1,n,1}(\mathbf{y})\sim\cdots\sim\mathcal{V}_{1,n,n-3}(\mathbf{y})$} edge [-] (a);
\node(c)[right = of a]{\ \ \ };
\node(d)[left = of a]{\ \ \ };
    \end{tikzpicture}}
\end{figure}\normalsize

\item[(\romannumeral2)] Assume $a=2$. Then $|y_1|=y_2=\cdots=y_{n-2}=1\leq y_{n-1}\leq y_n$, and for $j\leq n-2$,
\begin{align*}
        |2y_1\cdots y_{j-1}y_{j+1}\cdots y_n-y_j|&\geq2|y_{n-1}y_n|-1\geq|y_j|,\\
|2y_1\cdots y_{n-2}y_n-y_{n-1}|&\geq2|y_n|-|y_{n-1}|\geq|y_{n-1}|,
    \end{align*}
Note that $|2y_1\cdots y_{j-1}x_{j+1}\cdots y_n-y_j|=|y_j|$ implies $y_1=\cdots=y_n=1$, and that $2y_1y_n-y_{n-1}=y_{n-1}$ implies $y_1=\cdots=y_{n-2}=1\leq y_{n-1}=y_n$. Thus, there are four cases, as shown in the following diagrams:\footnotesize
\begin{figure}[H]
    \centering
\subcaptionbox{\footnotesize$\Delta\mathcal{V}_{2,n,n}(\mathbf{y})<\Delta(\mathbf{y})<\Delta\mathcal{V}_{2,n,n-1}(\mathbf{y})$\normalsize}{\begin{tikzpicture}[every node/.style={align=center}] 
    \node(a){$\mathbf{y}$}; 
    \node(b)[below =of a]{$\mathcal{V}_{2,n,n}(\mathbf{y})$} edge [-](a);
    \node(c)[above= of a]{$\cdots$} edge [-] (a);
    \node(d)[above left= of a]{$\mathcal{V}_{2,n,1}(\mathbf{y})$} edge [-] (a);
    \node(e)[above right= of a]{$\mathcal{V}_{2,n,n-1}(\mathbf{y})$} edge [-] (a);
    \end{tikzpicture}} 
\hspace{8pt} 
\subcaptionbox{\footnotesize$\Delta(\mathbf{y})<\Delta\mathcal{V}_{2,n,n-1}(\mathbf{y}),\Delta\mathcal{V}_{2,n,n}(\mathbf{y})$\normalsize}{\begin{tikzpicture}[every node/.style={align=center}] 
    \node(a){$\mathbf{y}$};
    \node(c)[above= of a]{$\cdots$} edge [-] (a);
    \node(d)[above left= of a]{$\mathcal{V}_{2,n,1}(\mathbf{y})$} edge [-] (a);
    \node(e)[above right= of a]{$\mathcal{V}_{2,n,n}(\mathbf{y})$} edge [-] (a);
    \end{tikzpicture}}
\vspace{10pt}

\subcaptionbox{\footnotesize$\Delta(\mathbf{y})=\Delta\mathcal{V}_{2,n,n-1}(\mathbf{y})<\Delta\mathcal{V}_{2,n,n-2}(\mathbf{y})$\normalsize}{\begin{tikzpicture}[every node/.style={align=center}] 
    \node(a){$\mathbf{y}=(1,\dots,1,x_{n-1},x_{n-1})$}; 
    \node(b)[above= of a]{$\mathcal{V}_{2,n,1}(\mathbf{y})\sim\cdots\sim\mathcal{V}_{2,n,n-2}(\mathbf{y})$} edge [-] (a);
\node(c)[left= of a]{\ \ \ };
\node(d)[right= of a]{\ \ \ };
    \end{tikzpicture}}
\subcaptionbox{\footnotesize$\Delta(\mathbf{y})=\Delta\mathcal{V}_{2,n,n-1}(\mathbf{y})=\Delta\mathcal{V}_{2,n,n-2}(\mathbf{y})$\normalsize}{\begin{tikzpicture}[every node/.style={align=center}] 
    \node(a){$\mathbf{y}=(1,\dots,1)$};
\node(c)[left= of a]{\ \ \ \ \ \ \ };
\node(d)[right= of a]{\ \ \ \ \ \ \ };
    \end{tikzpicture}}
\end{figure}\normalsize
\end{enumerate}
\end{enumerate}
By \romannumeral1) $\sim$ \romannumeral3), we conclude that $\mathbf{y}\in\cup\mathfrak{T}_0'(a,k,n)\cup\mathfrak{T}_1'(a,k,n)\cup\mathfrak{T}_2'(a,k,n)\cup\mathfrak{T}_{>2}'(a,k,n)$ is either \begin{enumerate}
        \item[(\romannumeral1)] a descendent vertex of $\mathbf{y},\mathcal{V}_{a,n,i}(\mathbf{y})$ for all $i=1,\dots,n$, or 
        \item[(\romannumeral2)] the last vertex of the graph associated with $\Delta$, corresponding to the orbit $\Gamma_{a,n}\cdot\mathbf{x}$,
    \end{enumerate} 
which means that $\Gamma_{a,n}\cdot\mathbf{x}$ has exactly one last vertex. (If $\Gamma_{a,n}\cdot\mathbf{x}$ has last vertices $\mathbf{z}_1\neq\mathbf{z}_2$, then the orbit must have a vertex $\mathbf{z}$ such that $\mathbf{z}_1$ and $\mathbf{z}_2$ are descendent vertices of $\mathbf{z}$. By \romannumeral1) $\sim$ \romannumeral3), it follows that $a=1$ and $\mathbf{z}\in\mathfrak{T}_1'(1,k,n)$. However, Lemma \ref{2.2} implies that $\mathbf{z}_1=\mathbf{z}_1$, which is a contradiction.)\qed

 \textit{Proof of the main theorem}. Now, we will collect all of the last vertices of each partition $\mathfrak{T}_0'(a,k,n)$, $\mathfrak{T}_1'(a,k,n)$, $\mathfrak{T}_2'(a,k,n)$, and $\mathfrak{T}_{>2}'(a,k,n)$ of $\{\mathbf{x}\in V_{a,k,n}(\mathbb{Z}):|x_1|\leq x_2\leq\cdots\leq x_n\}$. Then all elements in this set of the last vertices are not $\Gamma_{a,n}$-equivalent by Lemma \ref{2.3}.
\begin{enumerate}
    \item[\romannumeral1)] Assume $\mathbf{x}\in \mathfrak{T}_0'(a,k,n)$. That is, $|ax_1\cdots x_{n-2}|=0$ and $x_1=0$. Then we have
    \begin{align*}
        x_2^2+\cdots+x_n^2=k\ \text{and}\ \mathcal{V}_{a,n,i}:x_i\mapsto -x_i\ \text{for all}\ i.
    \end{align*}
    Hence, all elements in 
    \begin{align*}
        \mathfrak{S}_0(a,k,n)&=\{\mathbf{x}\in V_{a,k,n}(\mathbb{Z}):0=x_1\leq x_2\leq\cdots\leq x_n\}\\
        &=\{(0,x_2,\dots,x_n)\in\mathbb{Z}^n:0\leq x_2\leq\cdots\leq x_n, x_2^2+\cdots+x_n^2=k\}
    \end{align*}
    are inequivalent.
    
    \item[\romannumeral2)] Assume $\mathbf{x}\in \mathfrak{T}_1'(a,k,n)$. Then $|x_1|=x_2=\cdots=x_{n-2}=1\leq x_{n-1}\leq x_n$. For all $a>1$, $\mathfrak{S}_1(a,k,n)\subset\mathfrak{T}_1'(a,k,n)=\emptyset$, and by Lemma \ref{2.2} and the above remark, all elements in 
    \begin{align*}
        \mathfrak{S}_1(1,k,n)&=\{\mathbf{x}\in V_{1,k,n}(\mathbb{Z}):-x_1=x_2=\cdots=x_{n-2}=1\leq x_{n-1}\leq x_n\}\\
        &=\{(-1,1,\dots,1,x_{n-1},x_n)\in \mathbb{Z}^n:1\leq x_{n-1}\leq x_n, x_{n-1}^2+x_n^2+x_{n-1}x_n=k-n+2\}
    \end{align*}
    are inequivalent.
    
    \item[\romannumeral3)] Assume $\mathbf{x}\in \mathfrak{T}_2'(a,k,n)$. For all $a>2$, $\mathfrak{S}_2(a,k,n)\subset\mathfrak{T}_2'(a,k,n)=\emptyset$.
\begin{enumerate}
\item[(\romannumeral1)] If $a=1$, $k\neq n+1$, $x_1>0$, and $|x_n|\leq|x_1\cdots x_{n-1}-x_n|=|2x_{n-1}-x_n|$, then
    \begin{align*}
        x_n\leq2x_{n-1}-x_n\ &\Rightarrow\ x_n\leq x_{n-1},\\
        x_n\leq-2x_{n-1}+x_n\ &\Rightarrow\ x_{n-1}\leq0,
    \end{align*}
    implying that $x_{n-1}=x_n$. If $x_1<0$, then $|x_1\cdots x_{n-1}-x_n|=|-2x_{n-1}-x_n|>|x_n|$. Consequently, by Lemma \ref{2.3}, all elements in 
    \begin{align*}
        \mathfrak{S}_2(1,k,n)&=\{\mathbf{x}\in V_{1,k,n}(\mathbb{Z}):x_1=x_2=\cdots=x_{n-3}=1<x_{n-2}=2\leq x_{n-1}= x_n\}\\
        &\ \ \ \ \cup\{\mathbf{x}\in V_{1,k,n}(\mathbb{Z}):-x_1=x_2=\cdots=x_{n-3}=1<x_{n-2}=2\leq x_{n-1}\leq x_n\}\\
        &=\{(1,\dots,1,2,x_{n-1},x_{n-1})\in\mathbb{Z}^n:2\leq x_{n-1},k=n+1\}\\
        &\ \ \ \ \cup\{(-1,1,\dots,1,2,x_{n-1},x_n)\in \mathbb{Z}^n:2\leq x_{n-1}\leq x_n, x_{n-1}+x_n=\sqrt{k-n-1}\}
    \end{align*}
    are inequivalent.
\item[(\romannumeral2)] If $a=2$, $k\neq n-2$, $x_1>0$, and $|x_n|\leq|2x_1\cdots x_{n-1}-x_n|=|2x_{n-1}-x_n|$, then
    \begin{align*}
        x_n\leq2x_{n-1}-x_n\ &\Rightarrow\ x_n\leq x_{n-1},\\
        x_n\leq-2x_{n-1}+x_n\ &\Rightarrow\ x_{n-1}\leq0,
    \end{align*}
    implying $x_{n-1}=x_n$. If $x_1<0$, then $|2x_1\cdots x_{n-1}-x_n|=|-2x_{n-1}-x_n|>|x_n|$. Consequently, by Lemma \ref{2.3}, all elements in 
    \begin{align*}
        \mathfrak{S}_2(2,k,n)&=\{\mathbf{x}\in V_{2,k,n}(\mathbb{Z}):x_1=x_2=\cdots=x_{n-2}=1\leq x_{n-1}= x_n\}\\
        &\ \ \ \ \cup\{\mathbf{x}\in V_{2,k,n}(\mathbb{Z}):-x_1=x_2=\cdots=x_{n-2}=1\leq x_{n-1}\leq x_n\}\\
        &=\{(1,\dots,1,x_{n-1},x_{n-1})\in\mathbb{Z}^n:1\leq x_{n-1},k=n-2\}\\
        &\ \ \ \ \cup\{(-1,1,\dots,1,x_{n-1},x_n)\in \mathbb{Z}^n:1\leq x_{n-1}\leq x_n, x_{n-1}+x_n=\sqrt{k-n+2}\}
    \end{align*}
    are inequivalent.
\end{enumerate}
   
    \item[\romannumeral4)] Assume $\mathbf{x}\in\mathfrak{T}_{>2}'(a,k,n)\subset V_{a,k,n}(\mathbb{Z})\setminus\mathfrak{T}(a,k,n)$.
    \begin{enumerate}
        \item[(\romannumeral1)] If $x_1>0$ and $x_n=|x_n|\leq|ax_1\cdots x_{n-1}-x_n|$, it follows that $x_n\leq\frac{a}{2}x_1\cdots x_{n-1}$.
        \item[(\romannumeral2)] If $x_1<0$, then $|x_1\cdots x_{n-1}-x_n|=x_1\cdots x_{n-1}+x_n>x_n=|x_n|$
    \end{enumerate}
    By (\romannumeral1), (\romannumeral2), and Lemma \ref{2.1}, all elements in
    \begin{align*}
        \mathfrak{S}_{>2}(a,k,n)&=\left\{\mathbf{x}\in V_{a,k,n}(\mathbb{Z}):0<x_1\leq\cdots\leq x_n\leq\frac{a}{2}x_1\cdots x_{n-1}\ \text{and}\ ax_1\cdots x_{n-2}>2\right\}\\
        &\ \ \ \ \cup \{\mathbf{x}\in V_{a,k,n}(\mathbb{Z}):0<-x_1\leq x_2\leq \cdots\leq x_n\ \text{and}\ ax_1\cdots x_{n-2}<-2\}
    \end{align*}
    are inequivalent.

\end{enumerate}
Therefore, we conclude that the set
\begin{align*}
    \mathfrak{S}_0(a,k,n)\cup\mathfrak{S}_1(a,k,n)\cup\mathfrak{S}_2(a,k,n)\cup\mathfrak{S}_{>2}(a,k,n)
\end{align*}
is a fundamental domain for $\Gamma_{a,n}$.\qed
\begin{corollary} If $(a,k)\neq (1,n+1),(2,n-2)$, then $V_{a,k,n}(\mathbb{Z})$ has finitely many $\Gamma_{a,n}$-orbits. Otherwise, it has infinitely many $\Gamma_{a,n}$-orbits.
\end{corollary} 
\textit{Proof}. Let $a\in\mathbb{Z}^+$ and $k\in\mathbb{Z}$ such that $(a,k)\neq (1,n+1),(2,n-2)$. We need to show that $\mathfrak{S}_0(k,n)\cup\mathfrak{S}_1(k,n)\cup\mathfrak{S}_2(k,n)\cup\mathfrak{S}_{>2}(k,n)$ is finite. Clearly, $\mathfrak{S}_0(a,k,n)$ and $\mathfrak{S}_1(a,k,n)$ are finite. Since $(a,k)\neq (1,n+1),(2,n-2)$, it follows that
\begin{align*}
    &\{(1,\dots,1,2,x_{n-1},x_{n-1})\in\mathbb{Z}^n:2\leq x_{n-1},k=n+1\}=\emptyset\ \text{if}\ k\neq n+1,\\
 &\{(1,\dots,1,x_{n-1},x_{n-1})\in\mathbb{Z}^n:1\leq x_{n-1},k=n-2\}=\emptyset\ \text{if}\ k\neq n-2,\\
&\mathfrak{S}_2(1,k,n)=\{(-1,1,\dots,1,2,x_{n-1},x_n)\in \mathbb{Z}^n:2\leq x_{n-1}\leq x_n, x_{n-1}+x_n=\sqrt{k-n-1}\},\\
    &\mathfrak{S}_2(2,k,n)=\{(-1,1,\dots,1,x_{n-1},x_n)\in \mathbb{Z}^n:1\leq x_{n-1}\leq x_n, x_{n-1}+x_n=\sqrt{k-n+2}\}.
\end{align*}
Thus, $\mathfrak{S}_2(1,k,n)$ and $\mathfrak{S}_2(2,k,n)$ are finite. Now, choose $\mathbf{x}\in\mathfrak{S}_{>2}(a,k,n)$. 
\begin{enumerate}
    \item[(\romannumeral1)] If $x_1<0$, then
\begin{align*}
    k=x_1^2+\cdots+x_n^2+|ax_1\cdots x_n|\geq x_n^2.
\end{align*}
Thus, we obtain $0<-x_1\leq \cdots\leq x_n\leq\sqrt{k}$. 
    \item[(\romannumeral2)] Assume $x_1>0$. Note that $0<x_i\leq\cdots\leq x_{n-1}$ are bounded. (If there is  $1\leq i\leq n-2$ such that $M\leq x_i\leq\cdots\leq x_n$ for all $M>0$, since $x_1\cdots x_{n-1}\geq Mx_{n-1}$, then we have
\begin{align*}
k&=x_1^2+\cdots+x_n^2-ax_1\cdots x_n\\
&\leq (n-1)x_{n-1}^2+x_n^2-aMx_{n-1}x_n\\
&\approx \left(\sqrt{n-1}x_{n-1}-\frac{aMx_n}{2\sqrt{n-1}}\right)^2-\left(\frac{aMx_n}{2\sqrt{n-1}}\right)^2\\
&=x_{n-1}((n-1)x_{n-1}-aMx_n)<k
\end{align*} 
for some sufficiently large $M>n$, which is a contradiction.) Say $x_1\leq\cdots\leq x_{n-2}\leq M$. When $k-x_1^2+\cdots+x_{n-2}\neq0$, since  $ax_1\cdots x_{n-2}>2$, then the equation
    \begin{align*}
        \left(x_n-\frac{a}{2}x_1\cdots x_{n-1}\right)^2-\left(\left(\frac{a}{2}x_1\cdots x_{n-2}\right)^2-1\right)x_{n-1}^2=k-x_1^2-\cdots-x_{n-2}^2
    \end{align*}
    is shown as follows: 
    \begin{figure}[H]
\centering
\subcaptionbox{\footnotesize$k-x_1^2-\cdots-x_{n-2}^2>0$\normalsize}{\includegraphics[width=0.492\textwidth]{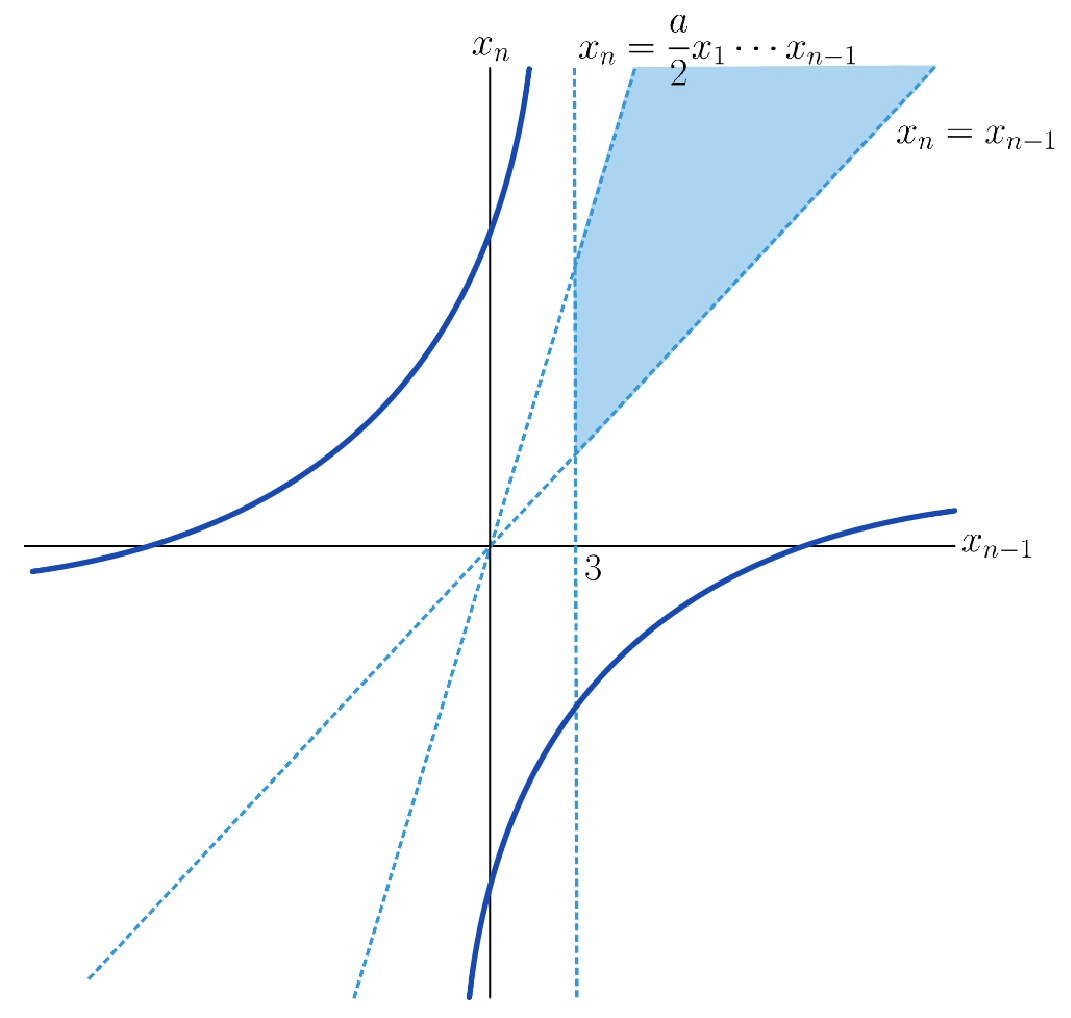}}%
\hspace{3pt} 
\subcaptionbox{\footnotesize$k-x_1^2-\cdots-x_{n-2}<0$\normalsize}{\includegraphics[width=0.492\textwidth]{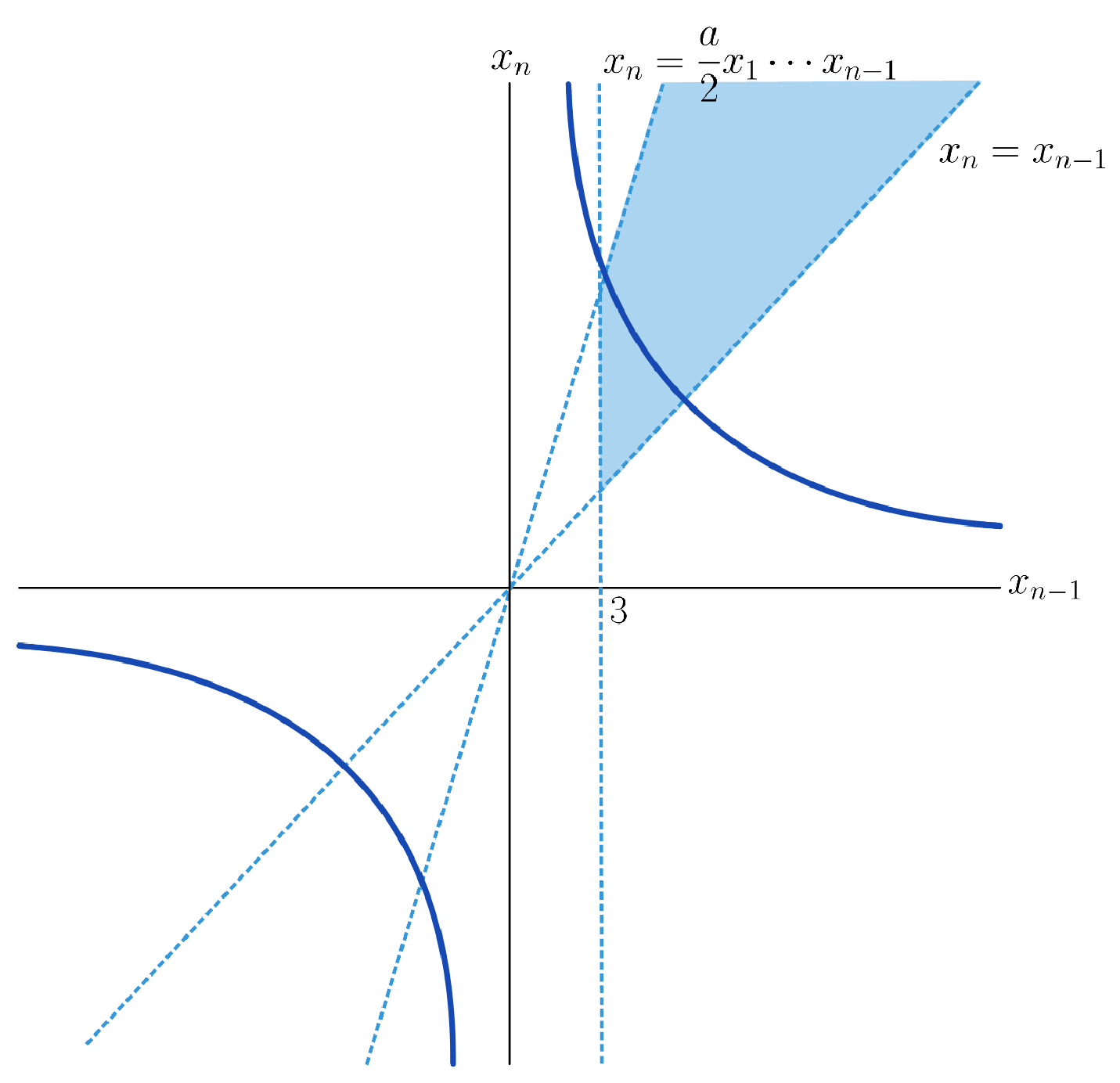}}%
\caption{\small Curves when $k-x_1^2-\cdots-x_{n-2}^2\neq0$\normalsize}
\end{figure}
\noindent As shown in Figure 5, when $k-x_1^2-\cdots-x_{n-2}^2>0$, since  
\begin{align*}
-\left(\left(\frac{a}{2}x_1\cdots x_{n-2}\right)^2-1\right)x_{n-1}^2=k-x_1^2-\cdots-x_{n-2}^2<0\ \ &\text{if}\ x_n=\frac{a}{2}x_1\cdots x_{n-1},\\
2x_{n-1}^2-ax_1\cdots x_{n-2}x_{n-1}^2=k-x_1^2-\cdots-x_{n-2}^2<0\ \ &\text{if}\ x_n=x_{n-1},
\end{align*}
it follows that the curve cannot intersect either $x_n=\frac{a}{2}x_1\cdots x_{n-1}$ or  $x_n=x_{n-1}$. That is, $\mathbf{x}$ does not belong to $\mathfrak{S}_{>2}(a,k,n)$. When $k-x_1^2-\cdots-x_{n-2}^2<0$, since 
the curve intersects $x_n=\frac{a}{2}x_1\cdots x_{n-1}$ and  $x_n=x_{n-1}$ at $\left(2\sqrt{\frac{x_1^2+\cdots+x_{n-2}^2-k}{(ax_1\cdots x_{n-2})^2-4}}, ax_1\cdots x_{n-2}\sqrt{\frac{x_1^2+\cdots+x_{n-2}^2-k}{(ax_1\cdots x_{n-2})^2-4}}\right)$ and $\left(\sqrt{\frac{x_1^2+\cdots+x_{n-2}^2-k}{ax_1\cdots x_{n-2}-2}},\sqrt{\frac{x_1^2+\cdots+x_{n-2}^2-k}{ax_1\cdots x_{n-2}-2}}\right)$, respectively, there are only finitely many points belonging to $\mathfrak{S}_{>2}(a,k,n)$. If $k-x_1^2-\cdots-x_{n-2}^2=0$, then we obtain
\begin{align*}
    &0\leq\frac{a}{2}x_1\cdots x_{n-1}-x_n=\sqrt{\left(\frac{a}{2}x_1\cdots x_{n-2}\right)^2-1}x_{n-1}\\
    &\Rightarrow\ \left(\frac{a}{2}x_1\cdots x_{n-2}-\sqrt{\left(\frac{a}{2}x_1\cdots x_{n-2}\right)^2-1}\right)x_{n-1}=x_n\geq x_{n-1}\\
    &\Rightarrow\ \frac{a}{2}x_1\cdots x_{n-2}-\sqrt{\left(\frac{a}{2}x_1\cdots x_{n-2}\right)^2-1}\geq1\\
    &\Rightarrow\ ax_1\cdots x_{n-2}\leq2,
\end{align*}
contradicting that $ax_1\cdots x_{n-2}>2$.
\end{enumerate}
By (\romannumeral1) and (\romannumeral2), $\mathfrak{S}_{>2}(a,k,n)$ is finite. If $(a,k)=(1,n+1)$ or $(a,k)=(2,n-2)$, since
\begin{align*}
\{(1,\dots,1,2,x_{n-1},x_{n-1})\in\mathbb{Z}^n:2\leq x_{n-1}\}&\subset\mathfrak{S}_2(1,n+1,n),\\
\{(1,\dots,1,x_{n-1},x_{n-1})\in\mathbb{Z}^n:1\leq x_{n-1}\}&\subset\mathfrak{S}_2(2,n-2,n),
\end{align*}
it follows that $V_{a,k,n}(\mathbb{Z})$ has infinitely many $\Gamma_{a,n}$-orbits. \qed

\subsection{Compatibility with the result for the Markoff equation} When $a=1$ and $n=3$, Theorem \ref{1.1} says that the set
     \begin{align*}
    \mathfrak{S}_0(1,k,3)\cup\mathfrak{S}_1(1,k,3)\cup\mathfrak{S}_2(1,k,3)\cup\mathfrak{S}_{>2}(1,k,3),
\end{align*}
where 
\begin{align*}
    \mathfrak{S}_0(1,k,3)&=\{(0,x_2,x_3)\in \mathbb{Z}^3:0\leq x_2\leq x_3, x_2^2+x_3^2=k\},\\
    \mathfrak{S}_1(1,k,3)&=\{(-1,x_2,x_3)\in \mathbb{Z}^3:1\leq x_2\leq x_3, x_2^2+x_3^2+x_2x_3=k-1
    \}\\
&=\{(-1,x_2,x_3)\in \mathbb{Z}^3:1\leq x_2\leq x_3, (2x_2+x_3)^2+3x_3^2=4(k-1)\},\\
    \mathfrak{S}_2(1,k,3)&=\{(2,x_2,x_2)\in \mathbb{Z}^3:2\leq x_2, k=4\}\cup\{(-2,x_2,x_3)\in\mathbb{Z}^3:2\leq x_2\leq x_3,x_2+x_3=\sqrt{k-4}\}\\
&=\{(2,x_2,x_2)\in \mathbb{Z}^3:2\leq x_2, k=4\}\cup\{(-2,x_2,x_3)\in\mathbb{Z}^3:2\leq x_2\leq x_3,(x_2+x_3)^2+4=k\},\\
    \mathfrak{S}_{>2}(1,k,3)&=\left\{(x_1,x_2,x_3)\in \mathbb{Z}^3:3\leq x_1\leq x_2\leq x_3\leq\frac{1}{2}x_1x_2, x_1^2+x_2^2+x_3^2-x_1x_2x_3=k\right\}\\
        &\ \ \ \ \cup \{(-x_1,x_2,x_3)\in \mathbb{Z}^3:3\leq x_1\leq x_2\leq  x_3, x_1^2+x_2^2+x_3^2+x_1x_2x_3=k\},
\end{align*}
is a fundamental domain for $\Gamma_{1,3}$ on $V_{1,k,3}(\mathbb{Z})$. If $k\geq5$ is generic (i.e.,\ $k$ is in none of the forms $k=u^2+v^2$, $4(k-1)=u^2+3v^2$, or $k=4+u^2$) or $k<0$, then $\mathfrak{S}_s(1,k,3)=\emptyset$ for all $s=0,1,2$.  Hence, Theorem \ref{1.1} for $(a,n)=(1,3)$ is compatible with Ghosh and Sarnak's result\cite{gs}.


\end{document}